\documentclass[preprint,12pt]{elsarticle}

\usepackage{amsmath,amssymb,amsfonts,theorem,color,makeidx,latexsym,epsfig,,subfigure}
\journal{Applied and Computational Harmonic Analysis}

\newtheorem{defn}{Definition}[section]
\newtheorem{lemma}[defn]{Lemma}

{\theorembodyfont{\rmfamily}
\newtheorem{ex}[defn]{Example}}
\newtheorem{thm}[defn]{Theorem}
\newtheorem{prop}[defn]{Proposition}
\newtheorem{cor}[defn]{Corollary}

\numberwithin{equation}{section}

\newcommand{\h}{{\cal H}}

\newcommand{\ltr}{ L^2(\mathbb R) }

\newcommand{\mn}{\mathbb N}
\newcommand{\mr}{\mathbb R}

\newcommand{\mz}{\mathbb Z}
\newcommand{\mc}{\mathbb C}

\def\bp{{\noindent\bf Proof. \ }}
\def\ep{\hfill$\square$\par\bigskip}

\def\bqs{\begin{equation}}
\def\eqs{\tag*{$\square$}\end{equation}\par\bigskip}

\def\la{\langle}
\def\ra{\rangle}
\def\ga{\gamma}

\def\gtk{\left\{g_k\right\}_{k=1}^\infty}

\def\bop{\begin{op}\rm}
\def\eop{\end{op}}

\def\bee{\begin{eqnarray}}
\def\ene{\end{eqnarray}}
\def\bes{\begin{eqnarray*}}
\def\ens{\end{eqnarray*}}
\def\bei{\begin{itemize}}
\def\eni{\end{itemize}}
\def\bt{\begin{thm}}
\def\et{\end{thm}}
\def\bc{\begin{cor}}
\def\ec{\end{cor}}
\def\bpr{\begin{prop}}
\def\epr{\end{prop}}
\def\bl{\begin{lemma}}
\def\el{\end{lemma}}
\def\bd{\begin{defn}}
\def\ed{\end{defn}}
\def\bex{\begin{ex}}
\def\enx{\end{ex}}
\def\bfi{\begin{fig}}
\def\efi{\end{fig}}

\def\inr{\int_{-\infty}^\infty}

\def\beeq{\begin{equation}}
\def\eneq{\end{equation}}
\def\besq{\begin{equation*}}
\def\ensq{\end{equation*}}

\def\cas{{\cal S}(\mr)}
\def\ltrn{L^2(\mr^n)}

\def\ftki{\{f_k\}_{k\in I}}
\def\gtki{\{g_k\}_{k\in I}}
\def\suki{\sum_{k\in I}}

\def\bDiff{\mathcal{D}}

\def\graph{\mathrm{graph\,}}
\newcommand{\re}{\mathrm{Re\,}}
\newcommand{\im}{\mathrm{Im\,}}
\newcommand{\eps}{\varepsilon}
\newcommand{\card}{\mathrm{card\,}}
\newcommand{\bS}{\mathbb{S}}

\begin{document}

\begin{frontmatter}



\title{Directional time-frequency analysis via  continuous frames\tnoteref{t1}}
\tnotetext[t1]{Research partially supported by a Start-Up Grant for Collaboration between EuroTech Universities funded by the Presidential Office, Technische Univerist\"at M\"unchen, Germany}

\author{ Ole Christensen\fnref{dtu}}
\author{ Brigitte Forster \fnref{passau}}
\author{ Peter Massopust \fnref{tum,hmgu}}
\fntext[dtu]{Department of Applied Mathematics and Computer Science, Technical University of Denmark,
Building 303, 2800 Lyngby, Denmark, ochr@dtu.dk}
\fntext[passau]{Fakult\"at f\"ur Informatik und Mathematik, Universit\"at Passau,
Innstr. 33, 94032 Passau, Germany, brigitte.forster@uni-passau.de}
\fntext[tum]{Zentrum Mathematik, M6, Technische Universit\"at M\"unchen, Boltzmannstr. 3, 85747 Garching, Germany,
massopust@ma.tum.de}
\fntext[hmgu]{Helmholtz Zentrum M\"unchen,Ingolst\"adter Landstra{\ss}e 1, 8764 Neuherberg, Germany}

\begin{abstract}
Grafakos and Sansing \cite{GS} have shown how to obtain directionally sensitive time-frequency decompositions in $L^2(\mr^n)$ based on Gabor systems in $\ltr;$ the key tool is the ``ridge idea," which lifts a function of one variable to a function of several variables. We generalize their result by showing that similar results hold starting with general frames for  $L^2(\mr),$  both in the setting of discrete frames and continuous frames. This allows to apply the theory for several other classes of frames, e.g., wavelet frames and shift-invariant systems.  We will consider applications to the Meyer wavelet and complex B-splines.  In the special case of wavelet systems we show how to discretize the representations using $\epsilon$-nets.
\end{abstract}


\end{frontmatter}


\section{Introduction} \label{146c}

Expansions of functions or signals as superpositions of basic building blocks with desired properties is one of the main tools in signal analysis. The expansions can be either in terms of an integral, a discrete sum, or a combination of both.

Many real-world signals depend on more than one variable. Depending on the type of expansion one is interested in, there are various ways to obtain such expansions. If an orthonormal basis for $\ltr$ is given, one can obtain an orthonormal basis for $L^2(\mr^n)$ via a simple tensor product, but this is highly inefficient.  Some of the standard methods to obtain expansions in $\ltr,$ e.g., wavelet frames of Gabor frames, have similar versions in $\ltrn,$ but they might not be optimal in order to detect features or special properties of the signal at hand. Other expansions are born in $L^2(\mr^n),$ typically for $n=2,3,$ e.g., caplets  \cite{HR}, ridgelets \cite{C, CD} and shearlets \cite{GL}; all of these can be considered as higher-dimensional wavelet-type systems with additional structure.

A different approach (parallel to the ridgelet construction) for Gabor systems was proposed by Grafakos and Sansing \cite{GS}. Starting with Gabor systems in $\ltr$, they developed a directionally sensitive Gabor-type expansion in $\ltrn$
using ridge functions. Two approaches were discussed in \cite{GS}: a discrete one, based on Gabor frames for $\ltr,$ and a semi-discrete version based on continuous Gabor systems generated by two non-perpendicular functions.

In this paper, we extend the main results  in \cite{GS} in various ways.
First, we observe that the above mentioned non-orthogonality places \cite{GS} in the setting of continuous frames, originally developed by
Ali et al. \cite{AAG1} resp. by Kaiser \cite{Ka}. Using techniques from frame theory, we then prove that the results in \cite{GS} have parallel versions starting with general frames for $\ltr,$ both in the discrete and the continuous setting.
The results are applied to the Meyer wavelet and complex B-splines.  In the special case of wavelet systems we show how to discretize the representations using $\epsilon$-nets.

In the rest of the introduction, we will introduce some notation and state the necessary facts about ridge functions and (continuous) frames. Then, in Section \ref{257a} we present the generalizations of the results in \cite{GS}. Semi-discrete representations of functions in $L^1(\mr) \cap L^2(\mr)$ are investigated in Section 3, where we also
apply the results to the Meyer wavelet and to complex B-splines. In the final Section 4, we obtain fully discrete representations for wavelet-type systems on bounded domains, by replacing the integral over the unit sphere
by an appropriately chosen $\eps$-net.

Some remarks concerning the notation:
Since we deal with functions in $L^1(\mr)$ and lift them to functions in $L^2(\mr^n)$, we need to consider inner products and the Fourier transform on different spaces. In general, for
functions $f\in L^1(\mr^m), \ m\in \mn,$ we define the Fourier  transform by
\bes \widehat{f} (\ga) := \int_{\mr^m} f(x)e^{-2\pi ix\cdot \ga}\, d\ga, \ \ga \in \mr^m,\ens where $x \cdot \ga$ denotes the canonical inner product on $\mr^m.$
We extend the Fourier transform to a unitary operator on $L^2(\mr^m)$ in the usual way. The inverse Fourier transform of a function $f$ will be denoted by $f^\lor.$ Also, for functions $f,g: \mathbb{R}^m \to \mathbb{C}$, $m\in \mn$,
we use the notation
\bee\label{integral}
\langle f, g \rangle := \int_{\mathbb{R}^m} f(x) \overline{g(x)}\, dx
\ene
whenever the right hand side converges. The unit sphere in $\mr^m$ will be denoted by $\bS^{m-1},$ and the Schwartz space of rapidly decreasing functions on $\mr$ by $\cas.$

\subsection{Ridge functions and the Radon transform}\label{subsect1}

%\textbf{Ich habe $\bDiff$ neu definiert, weil wir sonst einen Syntax-Overload mit dem Dilatationsoperator haben, der auch einen reellen Index hat. Eventuell ist $\partial$ oder so etwas besser.}

Let us now introduce the ``ridge procedure" that lifts functions of one variable to functions of several variables. Ridge functions were originally introduced by Pinkus \cite{P}. Our starting point is to extend the ordinary differential operator on $\mr$ to certain non-differentiable functions. In fact,
given $\alpha>0,$ define the differential operator $\bDiff^\alpha,$ acting on functions $h\in \cas$, by
\bee \label{126a}  \bDiff^\alpha(h):= (\widehat{h}(\cdot) | \cdot|^\alpha )^\lor;\ene
this definition clearly also makes sense for a large class of non-differentiable functions.

In the entire note, we use the following terminology, which relates
functions (written with lower case letters), the corresponding ridge functions
(written with a subscript), the action of the differential operator on the given function (written with capital letters), and the associated ridge function (written with capital letters and a subscript).

\bd \label{126h} Consider any function $g\in \cas.$

\bei \item[(i)] For $u\in \bS^{n-1},$ define the ridge function $g_u$ on $\mr^n$ by
\bee \label{12g} g_u(x):= g(u\cdot x), \ x\in \mr^n.\ene
\item[(ii)] Let
\bee \label{126b} G(s):=\bDiff^{\frac{n-1}{2}} \,g(s), \ s\in \mr.\ene
\item[(iii)]
For $u\in \bS^{n-1},$ define the weighted ridge function $G_u$ by
\bee \label{126c} G_u(x):= G(u \cdot x), \ x\in \mr^n.
\ene
\eni \ed

%\begin{rem}\label{rem1.3}
% For the existence of integrals of the type \eqref{integral} that involve functions as defined in \eqref{126b}, it is not necessary that $g\in \cas$. It suffices that $f$ and $g$  belong to $L^{2}(\mr)$ or that additionally one of them belongs to a Sobolev space $H^r(\mr)$ with $r > \frac{n-1}{2}$, i.e., an appropriate subspace of $L^{2}(\mr)$. We will encounter such a situation later in Example \ref{ex3.3} when we consider complex B-splines.
% \end{rem}
%
%
Given $u\in \bS^{n-1}$  the {\it Radon transform} of a function $f\in {\cal S}(\mr^n)$
 (in the direction $u$) is defined by
\bee \label{126d} R_uf (s):= \int_{u\cdot x=s} f(x)\, dx, \ s\in \mr.
\ene
The  Radon transform can be extended to a bounded operator from $L^1(\mr^n)$ to $L^1(\mr),$ see, e.g., \cite[p. 16 ff.]{N}. We also note that the {\it Fourier slice theorem} relates the (one-dimensional) Fourier transform of the Radon transform of a function $f\in L^1(\mr^n)$ to the ($n$-dimensional) Fourier transform of $f$ by the formula
\bes \widehat{R_u(f)}(\eta)= \widehat{f}(\eta u), \, \eta \in \mr, u\in \bS^{n-1}.\ens

The following lemma shows a close relation between ridge functions and the Radon transform.

\bl \label{126e}
For $f\in L^1(\mr^n),$  $g\in \cas,$ and $u\in \bS^{n-1},$
%the following relation holds between ridge functions and the Radon transform:
\bee \label{126f} \la f, g_u\ra=  \la R_uf, g\ra.\ene \el

\bp \bes \la f, g_u\ra & = & \int_{\mr^n} f(x) \overline{g(u\cdot x)}\, dx
 =  \inr \left( \int_{u\cdot x=s}  f(x) \overline{g(u\cdot x)}\, dx
 \right) \, ds   \\ & = & \inr \left(\int_{u\cdot x=s}  f(x) dx
 \right) \, \overline{g(s)}\, ds = \la R_uf, g\ra.\ens
 \ep

\subsection{Continuous frames}
In this section we review some of the known results about general continuous frames, as well as their concrete manifestations within Gabor analysis and wavelet theory.

\bd \label{111b} Let $\h$ be a complex Hilbert space and $M$ a measure space
with a positive measure $\mu$. A continuous frame
\index{continuous frame} is a family of vectors $\{f_k\}_{k\in M}$
for which the following hold:

\vspace{.1in}\noindent (i)  For all $f\in \h$, the mapping
$k\mapsto \la f,f_k\ra$ is a measurable function on $M$.

\vspace{.1in}\noindent (ii) There exist constants $A,B>0$ such
that \bes A \ ||f||^2\le \int_{M}|\la f,f_k\ra|^2\, d\mu(k) \le B \
||f||^2, \ \forall f\in \h. \ens
The  continuous frame  $\{f_k\}_{k\in M}$ is tight if we can choose $A=B.$\ed

For every continuous frame, there exists at least one {\it dual continuous frame,} i.e., a continuous frame $\{g_k\}_{k\in M}$ such that each
$f\in \h$ has the representation \bee \label{146a}
f  =   \int_M \la f,f_k\ra g_k \,d\mu(k);\ene
the integral in \eqref{146a} should be interpreted in the weak sense, i.e., as
\bee \label{146b}
\la f, g\ra  =   \int_M \la f,f_k\ra \la g_k,g \ra \,d\mu(k), \ \forall f,g\in \h.\ene
If $\{f_k\}_{k\in M}$ is a continuous tight frame with bound
$A,$ then $\{A^{-1} f_k\}_{k\in M}$ is a dual continuous frame.

Continuous frames generalize the more widely known (discrete) frames. In fact, in the case where $I$ ($=M$) is a countable set equipped with the counting measure, Definition \ref{111b} yields the classical frames. Continuous frames were introduced independently by
Ali et al. \cite{AAG1} and Kaiser \cite{Ka}. We will present the most important concrete cases below; for constructions of (discrete) frames, see the monographs
\cite{daubechies,G,CBN}.

There are also several well known examples of continuous frames for $\ltr$ available in the literature. In order to introduce these, consider the translation, modulation, and scaling--operators on $\ltr$ defined by
\begin{eqnarray*}
T_af(x) := f(x-a), \ \
E_bf(x) := e^{2\pi i bx}f(x),\ \
D_c f(x) := c^{1/2}f(cx),
\end{eqnarray*} where $a,b\in \mr, c>0.$

A system of functions of the form $\{E_{b}T_ag\}_{a,b\in \mr}$ is called a
(continuous) {\it Gabor system.} We state the following well known result, see, e.g.,
\cite[Theorem 3.2.1]{G},\cite[Proposition 9.9.1]{CBN}.

\bpr \label{cgt5} Let $f_1,f_2,g_1,g_2\in \ltr$. Then
\bes
\inr \inr \la f_1, E_bT_ag_1\ra \overline{\la f_2, E_bT_ag_2\ra}\,db\,da
= \la f_1, f_2\ra \la g_2, g_1\ra.
\ens \epr
Proposition \ref{cgt5} has an immediate and well known consequence concerning the construction of continuous tight Gabor frames and dual pairs. The result shows that it is very easy to construct such frames, especially
 with windows belonging to the Schwartz space $\cas.$

\bc \label{cgt7} \

\bei \item[(i)] For any $g\in \ltr \setminus \{0\},$ the Gabor system
$\{E_bT_ag\}_{a,b\in \mr}$ is a continuous tight frame for $\ltr$ with
respect to $M= \mr^2$ equipped with the Lebesgue measure, with frame bound
$A=||g||^2.$
\item[(ii)] For any functions $g_1, g_2\in \ltr$ for which $\la g_1, g_2\ra \neq 0,$ the Gabor systems  $\{E_{b}T_ag_1\}_{a,b\in \mr}$ and
  $\{\frac1{\la g_1, g_2\ra}\, E_{b}T_ag_2\}_{a,b\in \mr}$  are dual continuous frames. \eni \ec

A {\it wavelet system} has the form   $\{D_aT_b\psi\}_{a\neq 0, b\in \mr}$
for a suitable function $\psi\in \ltr.$ We say that $\psi$ satisfies the {\it
admissibility condition} \index{admissibility condition} if \bee
\label{ad1} C_\psi:= \inr \frac{|\widehat{\psi}(\gamma)|^2}{|\gamma |}
d\gamma < \infty. \ene
The admissibility condition gives the following result, see, e.g., \cite[Prop. 2.4.1]{daubechies}:

\bpr \label{cwt7}
Assume that $\psi$ is admissible. Then, for all functions $f,g\in \ltr$,
\bee \label{cwt2}
\inr \inr \la f, D_aT_b\psi\ra \overline{ \la g, D_aT_b\psi\ra }\,\frac{da\, db}{a^2}
= C_\psi \la f, g\ra.
\ene

\epr

Again, Proposition \ref{cwt7} immediately leads to a construction of a tight frame:

\bc \label{cwt4} If $\psi\in \ltr$ is admissible, then
$\{D_aT_b \psi\}_{a\neq 0, b\in \mr}$ is a continuous frame for
$\ltr$ with respect to $\mr \times \mr \setminus \{0\}$ equipped
with the Haar measure $\frac{1}{a^2}da\,db,$  with frame bound
$A= C_\psi.$ \ec

Note that the admissibility condition is easy to satisfy, even with
generators $\psi \in \cas.$ In fact, all functions $\psi \in \mathcal{S}(\mathbb{R})$ with vanishing mean
$$
\int_{\mathbb{R}} \psi(x) \, dx = \widehat{\psi}(0) =0
$$
satisfy the admissibility condition, see, e.g., \cite{lmr94}.

\section{Decompositions via continuous frames}\label{257a}

We first show that any pair of continuous dual frames for $\ltr$ consisting of functions in $\cas$ leads to an integral representation of
functions $f\in L^1(\mr^n)$ for which
$\widehat{f} \in L^1(\mr^n).$  This generalizes Theorem 3 in \cite{GS}. Note that Theorem 3 in \cite{GS} does not use the term {\it continuous frame,}
but just the technical condition $\la g_1, g_2\ra\neq 0;$ in the particular context of Gabor analysis this means exactly that the functions $g_1, g_2$ generate
continuous dual Gabor frames, as we saw in Corollary \ref{cgt7}.

 \bt \label{146d} Let $\{f_k\}_{k\in M}$ and $\{g_k\}_{k\in M}$ be dual
continuous frames for $\ltr,$ consisting of functions in $\cas.$ Then, for $f\in L^1(\mr^n)$ such that
$\widehat{f} \in L^1(\mr^n),$
\bee \label{136a} f= \frac12 \int_{\bS^{n-1}} \int_M \la f, G_{k,u}\ra F_{k,u} \, dk \, du.\ene \et

\bp Calculating the left-hand-side using Lemma \ref{126e} yields
\bes & \ & \frac12 \int_{\bS^{n-1}} \int_M \la f, G_{k,u}\ra F_{k,u}(x) \, dk \, du \\ & = & \frac12 \int_{\bS^{n-1}} \int_M \la R_uf, G_{k}\ra F_{k}( u\cdot x) \, dk \, du \\ & = & \frac12 \int_{\bS^{n-1}} \int_M \la \widehat{R_uf}, \widehat{G_{k}}\ra F_{k}( u\cdot x) \, dk \, du \\ & = & \frac12 \int_{\bS^{n-1}} \int_M \inr \widehat{R_uf}(\sigma)\overline{ \widehat{g_{k}}(\sigma)}\, |\sigma|^{\frac{n-1}{2}}\, d\sigma F_{k}( u\cdot x) \, dk \, du
.\ens
Now,
\bes F_k(u\cdot x)= \left(\widehat{f_k}(\cdot) | \cdot |^{\frac{n-1}{2}} \right)^{\lor}(u\cdot x)= \inr e^{2\pi i \eta u\cdot x}\widehat{f_k}(\eta) |\eta|^{\frac{n-1}{2}}\,d\eta, \ens  so we arrive at
\bes & \ & \frac12 \int_{\bS^{n-1}} \int_M \la f, G_{k,u}\ra F_{k,u}(x) \, dk \, du \\ & = & \frac12 \int_{\bS^{n-1}} \int_M \inr \widehat{R_uf}(\sigma)\overline{ \widehat{g_{k}}(\sigma)}\, |\sigma|^{\frac{n-1}{2}}\, d\sigma \inr e^{2\pi i \eta u\cdot x}\widehat{f_k}(\eta) |\eta|^{\frac{n-1}{2}}\,d\eta \, dk \, du \\ & = &
\frac12 \int_{\bS^{n-1}} \inr \left(\int_M  \inr \, |\sigma|^{\frac{n-1}{2}} \widehat{R_uf}(\sigma)\overline{ \widehat{g_{k}}(\sigma)}\, d\sigma \widehat{f_k}(\eta) \,dk\right) e^{2\pi i \eta u\cdot x} |\eta|^{\frac{n-1}{2}}\,d\eta \, du.
\ens Note that because $\{f_k\}_{k\in M}$ and $\{g_k\}_{k\in M}$ are dual
continuous frames for $\ltr,$ also $\{\widehat{f_k}\}_{k\in M}$ and $\{\widehat{g_k}\}_{k\in M}$ are dual continuous frames. The term above in the parantheses is exactly the frame decomposition with respect to these frames of the function
$| \cdot |^{\frac{n-1}{2}} \widehat{R_uf}(\cdot),$ evaluated at the point $\eta;$ thus
\bes \left( \int_M\inr \, |\sigma|^{\frac{n-1}{2}} \widehat{R_uf}(\sigma)\overline{ \widehat{g_{k}}(\sigma)}\, d\sigma \widehat{f_k}(\eta) \,dk\right)=|\eta|^{\frac{n-1}{2}} \widehat{R_uf}(\eta).\ens 
Inserting this yields that
\bee  \frac12 \int_{\bS^{n-1}} \int_M &&\hspace*{-24pt}\la f, G_{k,u}\ra F_{k,u}(x) \, dk \, du\nonumber\\
&  = &\frac12 \int_{\bS^{n-1}} \inr  |\eta|^{\frac{n-1}{2}} \widehat{R_uf}(\eta) e^{2\pi i \eta u\cdot x} |\eta|^{\frac{n-1}{2}}\,d\eta \, du \nonumber \\ 
\label{111a} & = &\frac12 \int_{\bS^{n-1}} \inr  |\eta|^{n-1} \widehat{R_uf}(\eta)
 e^{2\pi i \eta u\cdot x} \,d\eta \, du.
 \ene
From here on, we can just refer to the  proof of Theorem 3 in \cite{GS} in order to complete the proof. In fact, by the Fourier slice theorem, $\widehat{R_uf}(\eta)
= \widehat{f}(\eta u);$ inserting this in \eqref{111a}, and splitting the integral over $\mr$ into integrals over $]-\infty,0]$ and $[0, \infty[,$ a few changes of variables yield that \eqref{111a} equals $\int_{\mr^n} \widehat{f}(y)e^{2\pi i x\cdot y}\, dy=f(x),$ as desired.\ep

We have already sen in Section \ref{146c} that it is easy to construct continuous tight wavelet frames for $\ltr$ that are generated by
functions $\psi\in \cas;$ thus, it is easy to give applications of
Theorem \ref{146d}. However, for the purpose of applications our ultimate goal is to provide discrete realizations of the theory, so we will not consider concrete cases here.

\section{Semi-discrete representations}

The integral representation in Theorem \ref{146d} involves integrals over as well the sphere $\bS^{n-1}$ as the set $M.$ Letting $M$ be a discrete set equipped with the counting measure, we can of course apply the result to discrete frames as well; in this case we obtain what we will call a semi-discrete representation of
functions $f\in L^1(\mr^n) \cap L^2(\mr^n),$ only involving an integral over $\bS^{n-1}$ and a sum over the discrete index set. Nevertheless, we will follow the approach by
Grafakos and Sensing, see Theorem 5 in \cite{GS}, where a semi-discrete representation is derived in the Gabor case, independently of the continuous case. The reason for doing this is that the technical conditions are slightly weaker in this approach, leading to a representation that is valid for a larger class of functions.

In order to prove the next theorem, we need a result that is stated as part of Lemma 2, \cite{GS}.

\bl
Given a function $f\in L^1(\mr^n)\cap L^2(\mr^2)$, we have that
\[
\bDiff^{\frac{n-1}{2}}R_uf \in L^2(\mr)
\]
for almost every $u\in \bS^{n-1}$.
\el

\bt \label{126k} Let $f\in L^1(\mr^n) \cap L^2(\mr^n)$ and let $I$ be a countable index set.
\bei \item[(i)] Let $\gtki \subset \cas$ denote a frame for $\ltr$ with frame bounds $A,B,$ and define the associated functions $G_k$ and $G_{k,u}$ as in Definition \ref{126h}. Then
\bee \label{126m} A\, ||f||^2 \le \int_{\bS^{n-1}} \suki | \la f, G_{k,u}\ra|^2 \, du \le B ||f||^2.\ene
\item[(ii)] Assuming that $\gtki$ and $\ftki$ are dual frames for $\ltr,$
both consisting of functions in $\cas,$ then
\bee \label{126n} f=  \frac12 \int_{\bS^{n-1}} \suki  \la f, G_{k,u}\ra F
_{k,u} \, du.\ene
\eni \et

\bp By Lemma \ref{126e}, applied to the function $G_k,$ we have that
\bee \label{126r} \la f, G_{k,u}\ra = \la R_uf, G_k\ra= \la R_uf, \bDiff^{\frac{n-1}{2}}g_k\ra = \la \bDiff^{\frac{n-1}{2}}R_uf, g_k\ra,\ene where the last equality following by partial integration and the assumption $g_k\in \cas.$
    Now, by the frame assumption on $\gtk,$
\bee \label{126q} A\, || \bDiff^{\frac{n-1}{2}}R_uf||^2 \le
\suki |\la  \bDiff^{\frac{n-1}{2}}R_uf, g_k\ra|^2  \le B\, || \bDiff^{\frac{n-1}{2}}R_uf||^2.\ene
As shown in \cite{GS}, \bee \label{126p} \int_{\bS^{n-1}}|| \bDiff^{\frac{n-1}{2}}R_uf||_{L^2(\mr)}^2 \, du= 2\, ||f||^2.\ene
Thus, integrating \eqref{126q} over $\bS^{n-1}$ and applying \eqref{126r} yields the result in (i).

As for the proof of (ii), the frame decomposition associated with the
frames $\ftki$ and $\gtki$ and applied to the function $\bDiff^{\frac{n-1}{2}}R_uf$ yields that
\bes \bDiff^{n-1}R_uf= \bDiff^{\frac{n-1}{2}}\bDiff^{\frac{n-1}{2}}R_uf
& = & \bDiff^{\frac{n-1}{2}} \suki \la \bDiff^{\frac{n-1}{2}}R_uf, g_k\ra f_k \\
& = &  \suki \la \bDiff^{\frac{n-1}{2}}R_uf, g_k\ra \bDiff^{\frac{n-1}{2}}f_k \\
& = &  \suki \la R_uf, G_k\ra F_k  \\
& = &  \suki \la f, G_{k,u}\ra F_k .\ens
Since
\bes f(x)=  \frac12 \int_{\bS^{n-1}} \bDiff^{n-1}R_u(f)(u\cdot x)\, du,\ens
it follows that
\bes f(x) & = & \frac12 \int_{\bS^{n-1}}   \suki \la f, G_{k,u}\ra F_k  (u\cdot x)\, du  = \frac12 \int_{\bS^{n-1}}  \suki \la f, G_{k,u}\ra  F_k  (u\cdot x)\, du \\ & = &  \frac12 \int_{\bS^{n-1}}  \suki \la f, G_{k,u}\ra  F_{k,u} (x)\, du,  \ens
as desired. \ep

It is easy to satisfy the assumptions in Theorem \ref{126k}; see, e.g., \cite[Theorem 3.4]{Wo}.
Let us illustrate the result by considering the Meyer wavelet.

\bex
Let $\nu: \mathbb{R} \to [\, 0,1\,]$ be a smooth function of sigmoidal shape required to satisfy
$ \nu(y) =0$ for $y \leq 0$, $\nu(y) =1$ for $y \geq 1$, and $ \nu(y) + \nu(1-y) =1$. An example of such a function is for instance the polynomial $\nu(y) = y^{4}(35-84y+70y^2-20 y^{3})$, for $y\in [\,0,1\,]$.

\begin{figure}[p]
\centering
\includegraphics[width = 0.45 \textwidth]{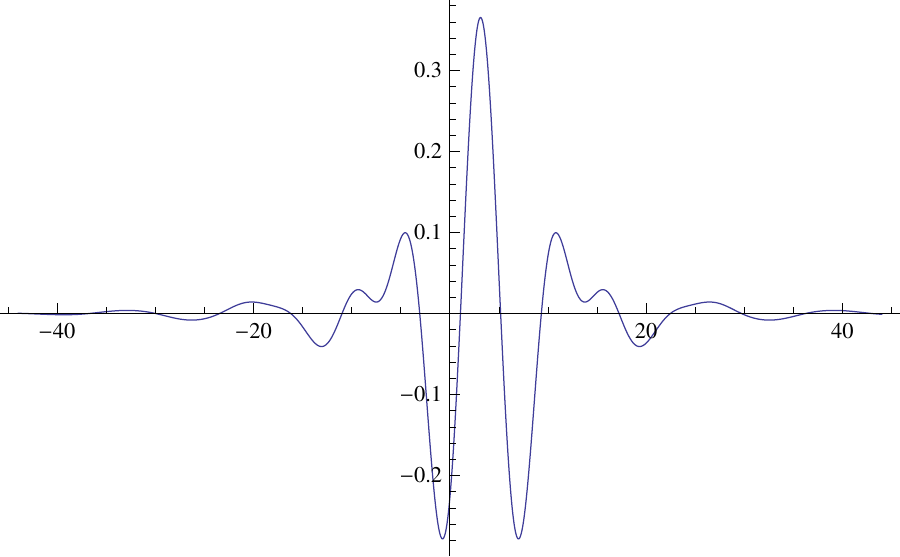}
\quad
\includegraphics[width = 0.45 \textwidth]{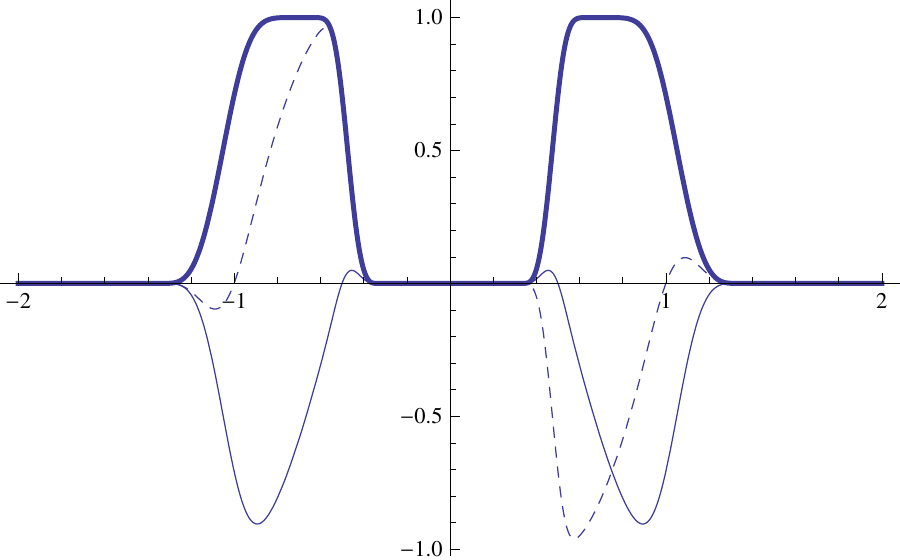}

\medskip
(a) \quad \quad (b)

\bigskip

\includegraphics[width = 0.45 \textwidth]{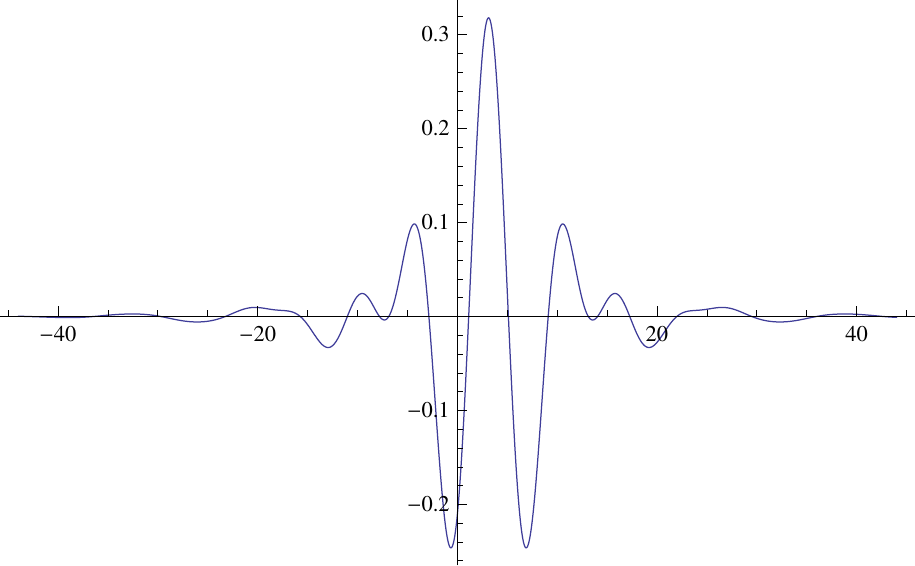}
\quad
\includegraphics[width = 0.45 \textwidth]{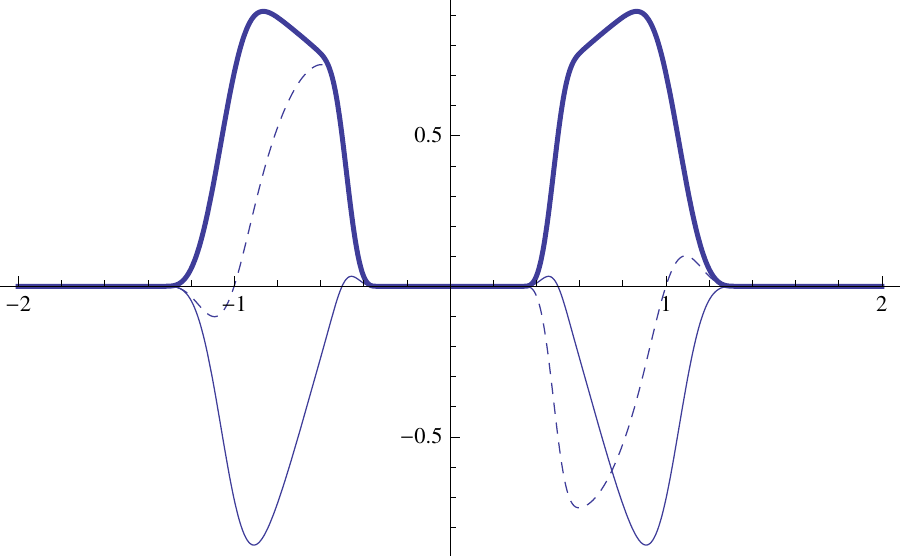}

\medskip

(c) \quad \quad (d)

\includegraphics[width = 0.7 \textwidth]{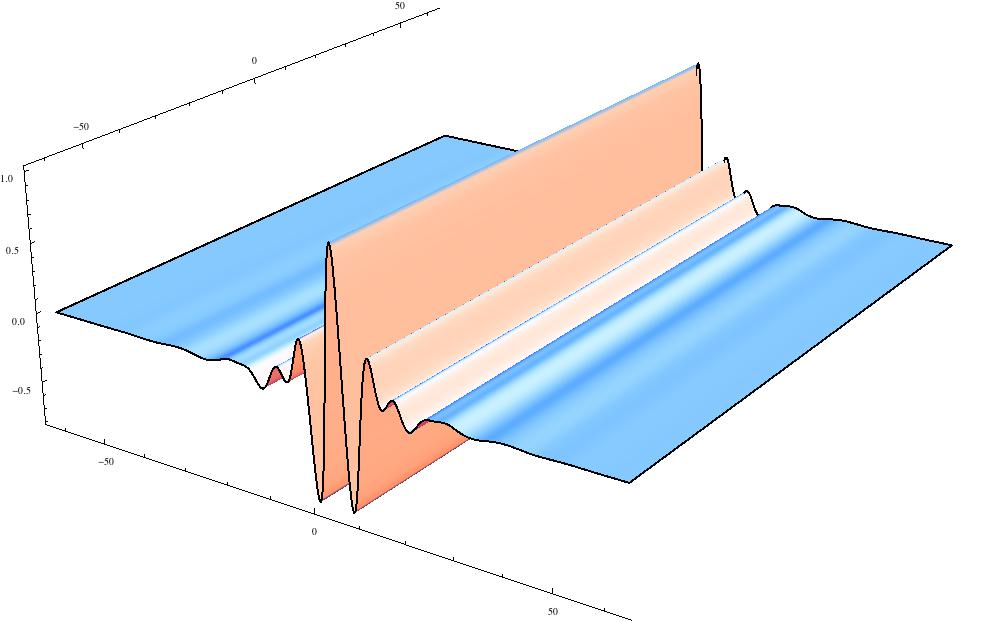}

(e)

\caption{(a) The real-valued Meyer wavelet $\psi_{0,0}$ and (b) its Fourier transform;  real, imaginary part and absolute value are plotted in normal, dashed and thick lines, respectively. (c) The generator $\Psi_{0,0}$ and (d) its Fourier transform (again real, imaginary part and absolute value). (e) The ridge frame generator $\Psi_{0,0,u}$ with $u = (1,0)^{T}$.}

\label{fig Meyer wavelet}
\end{figure}

Now let
$$
w(y) := \left\{
\begin{array}{c@{\quad}l}
\sin\left(\frac{\pi}{2}\nu (\frac{3y}{2\pi}-1)\right),& \mbox{for } \frac{2 \pi}{3} \leq y \leq \frac{4 \pi}{3},
\\[0.8ex]
\cos\left(\frac{\pi}{2}\nu (\frac{3y}{2\pi}-1)\right), & \mbox{for } \frac{4 \pi}{3} \leq y \leq \frac{8 \pi}{3},
\\[0.8ex]
0, & \mbox{elsewhere.}
\end{array}
\right.
$$
The classical Meyer wavelet $\psi$ is defined in the Fourier domain by
$$
\widehat{\psi}(\gamma) := e^{-i \pi \gamma} (w(2 \pi \gamma) + w(- 2 \pi \gamma)).
$$
It is well known that $\psi$ is a Schwartz function and that
$$
\{ \psi_{k,m} \}_{k,m \in \mathbb{Z}}:= \{2^{-m/2}\psi(2^{-m} \cdot -k) \mid k,m \in \mathbb{Z}\}
$$ is an orthonormal basis for $L^{2}(\mathbb{R}), $ see \cite{daubechies,lmr94,Wo}.
In particular,
$\{ \psi_{k,m} \}_{k,m \in \mathbb{Z}}$ is a frame, which is its own dual.
Thus, we can apply Theorem \ref{126k}; the functions $G_{k,u}=F_{k,u}$ have the form
\bee
\Psi_{k,m,u}(x) & := & \Psi_{k,m}(u \cdot x) = \bDiff^{\frac{n-1}{2}} \psi_{k,m}(u\cdot x) \nonumber\\
& = & (|\cdot|^{\frac{n-1}{2}}\widehat{\psi_{k,m}})^\lor (u \cdot x), \qquad k,m\in \mz, \; u\in\bS^{n-1}.\nonumber
\ene
For $n=2$ , the Meyer wavelet $\psi = \psi_{0,0}$ and the ridge function $\Psi_{0,0,u}$ with $u := (1,2)^{T}/\sqrt{5}$, are plotted in Figure \ref{fig Meyer wavelet}.
\enx

\bex\label{ex3.3}

%The motivation behind the definition of complex B-splines is the need for a single-band frequency analysis. For some applications, e.g., for phase retrieval tasks, complex-valued analysis bases are needed since real-valued bases can only provide a symmetric spectrum. Complex B-splines combine the advantages of spline approximation with an approximate one-sided frequency analysis, see equation (\ref{eq Frequenz-Produkt komplexe B-Splines}) \cite{forster06,fivegoodreasons}.

Complex B-splines $\beta_z: \mr \to\mc$ are a natural extension of the classical Curry-Schoenberg B-splines.  They were introduced in \cite{forster06} and then studied and extended in a series of papers, see, e.g.,  \cite{FGMS,FMUe,forstermasso10,multivarSplines}.
Given $z\in \mc$ with $\re z>1,$ the complex B-spline $\beta_z$ is defined in the Fourier domain by
\begin{equation}\label{Bz}
\widehat{\beta}_z(\gamma ) := \left( \frac{1-e^{-2\pi i\gamma}}{2\pi i\gamma}\right)^z.
\end{equation}
Setting $\Omega:\mr\to\mc$, $\gamma\mapsto\frac{1-e^{-2\pi i\gamma}}{2\pi i\gamma}$, one notices that
$$\graph\Omega \cap \{(0,y) \in\mr\times \mr \mid y<0\} = \emptyset,$$
implying that complex B-splines reside on the main branch of the complex logarithm and are thus well-defined.

\begin{figure}[p]
\centering

\includegraphics[width = 0.45 \textwidth]{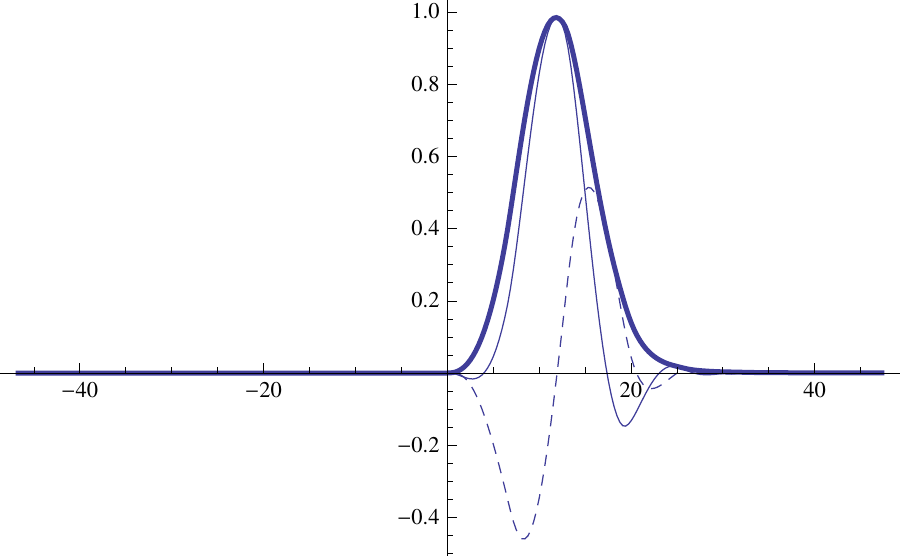}
\quad
\includegraphics[width = 0.45 \textwidth]{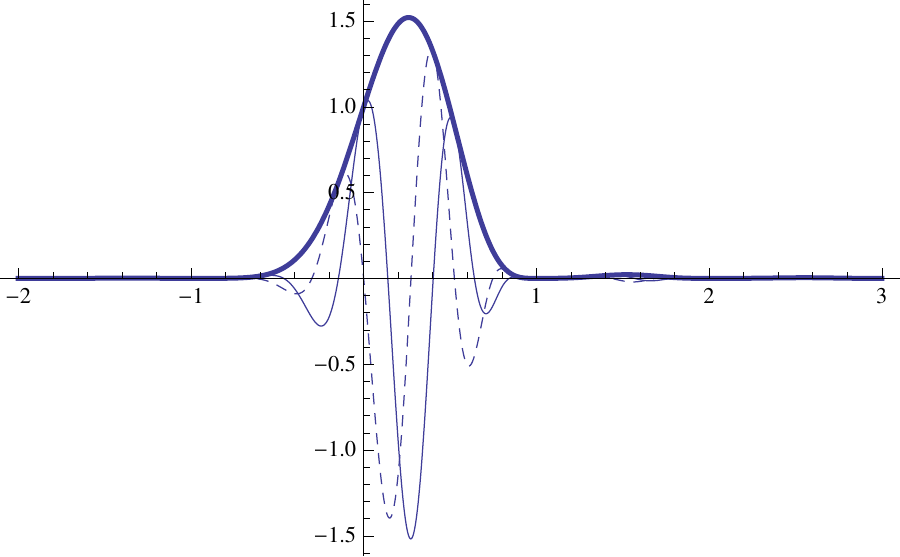}

\label{fig komplex Splines}
\caption{Complex B-spline $\beta_{z}$ in time domain (left) and frequency domain (right) for $z = 3.5 +i$. The spline in time domain has its support on the positive real axis, i.e., the spline is causal. Moreover, the frequency spectrum is shifted to the positive frequencies because of the positive imaginary part of $z$. Thick lines indicate the modulus, thin lines the real part and dashed lines the imaginary part.}
\end{figure}

Compared with the classical cardinal B-splines, complex B-splines $\beta_z$ possess an additional modulation and phase factor in the frequency domain:
\begin{equation}
\widehat{\beta}_z (\gamma) = \widehat{\beta}_{\re z}(\gamma)\,e^{i \im z \ln |\Omega(\gamma)|}\,e^{- \im z \arg \Omega(\gamma)}.
\label{eq Frequenz-Produkt komplexe B-Splines}
\end{equation}
The existence of these two factors allows the extraction of additional information from sampled data. In fact, the spectrum of a complex B-spline consists of the spectrum of a real-valued B-spline combined with a modulating and a damping factor. The presence of an imaginary part causes the frequency components on the negative and positive real axis to be  enhanced with different signs. This has the effect of shifting the frequency spectrum towards the negative or positive frequency side, depending on the sign of the imaginary part. This allows for an approximate single band analysis which is not possible with any real valued function, but necessary for certain phase retrieval tasks in signal processing\cite{forster06,fivegoodreasons}. Figure \ref{fig komplex Splines} gives an example of a complex B-spline in time domain as well as in frequency domain.

%The corresponding spline bases can be interpreted as approximate single-band filters \cite{forster06,fivegoodreasons}. In contrast to, for instance, complex wavelet bases, the phase information ($e^{i \im z \ln |\Omega(\gamma)|}$) is already built in, and an adjustable smoothness parameter, namely $\re z$, provides a continuous family of approximation and interpolation functions.

%The functions $\beta_z$ have, in general, non-compact support. However, the smoothness of their Fourier transform yields a fast decay in time domain:
%\[
%\beta_z (x) = \mathcal{O}(x^{-m}), \quad \mbox{for $\mn\ni m < \re z +1$,  as $x \to \infty$}.
%\]
%Furthermore, it was proved in \cite{forster06} that complex B-splines are elements of $L^1(\mr) \cap L^2(\mr).$  In fact,
For
$r< \re z - \frac12$ the complex B-splines
belong to the Sobolev spaces $H^r(\mr)$ (with respect to the $L^2$-Norm and with weight $(1+|x|^2)^r$).
%, due to their decay in frequency domain induced by the polynomial $\gamma^z$ in the denominator of \eqref{Bz}.
%It was also shown in \cite{forster06} that

Complex B-splines generate a multiresolution analysis $\{V_k \mid k\in \mz\}$ of $L^2(\mr)$.
In particular, $\{\beta_z (\cdot - \ell)\mid \ell\in \mz \}$ is a Riesz basis for $V_0$.

\begin{figure}[p]
\centering

\includegraphics[width = 0.45\textwidth]{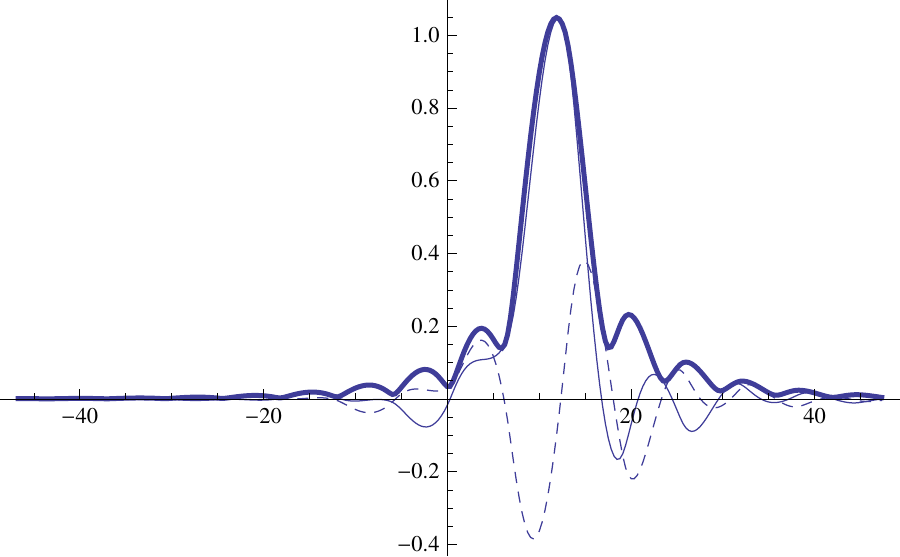}
\quad
\includegraphics[width = 0.45\textwidth]{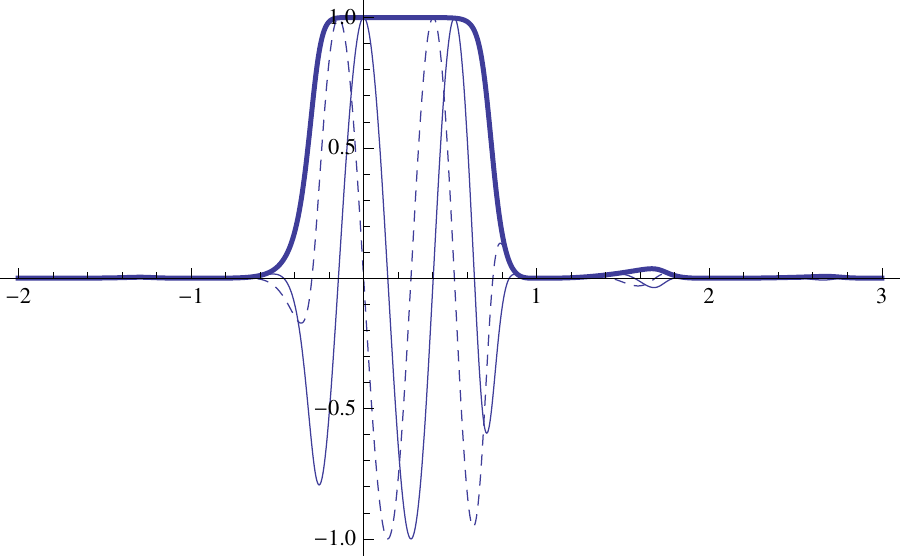}

\label{fig Spline ortho}
\caption{Complex orthonormalized B-spline $\beta_{z,\perp}$ in time domain (left) and frequency domain (right) for $z = 3.5 +i$. Thick lines indicate the modulus, thin lines the real part and dashed lines the imaginary part.
In comparison with figure \ref{fig komplex Splines}, $\beta_{z,\perp}$ is not causal anymore, but still has the shift in the frequency spectrum to the positive real axis.}
\end{figure}

To construct the corresponding orthonormalized wavelet $\psi_{z}$, we apply a high pass filter to the complex B-spline scaling function $\beta_{z}$. To this end, denote by
$$
A_{z}(\gamma) = \sum_{k\in\mathbb{Z}}|\widehat{\beta}_{z}(\gamma + k)|^{2}
$$
the autocorrelation filter. The convergence of the series is proved in \cite{forster06}.
Then
$$
\widehat{\beta}_{z,\perp}(\gamma) = \widehat{\beta}_{z}(\gamma) / \sqrt{A_{z}(\gamma)}
$$
is an orthonormal scaling function. For an example, see Figure \ref{fig Spline ortho}.

The scaling filter is denoted by
$$
H_{z}(\gamma/2) = \frac{\widehat{\beta}_{z,\perp}(\gamma)}{\widehat{\beta}_{z,\perp}(\gamma/2)}.
$$
The associated orthonormal wavelet $\psi_{z,\perp}$ is given by
$$
\widehat{\psi}_{z,\perp}(\gamma) = - e^{-i\pi\gamma} \, \overline{H_{z}((\gamma+1)/2} \,\, \widehat{\beta}_{z,\perp}(\gamma/2).
$$

\begin{figure}[p]
\centering

\includegraphics[width = 0.45 \textwidth]{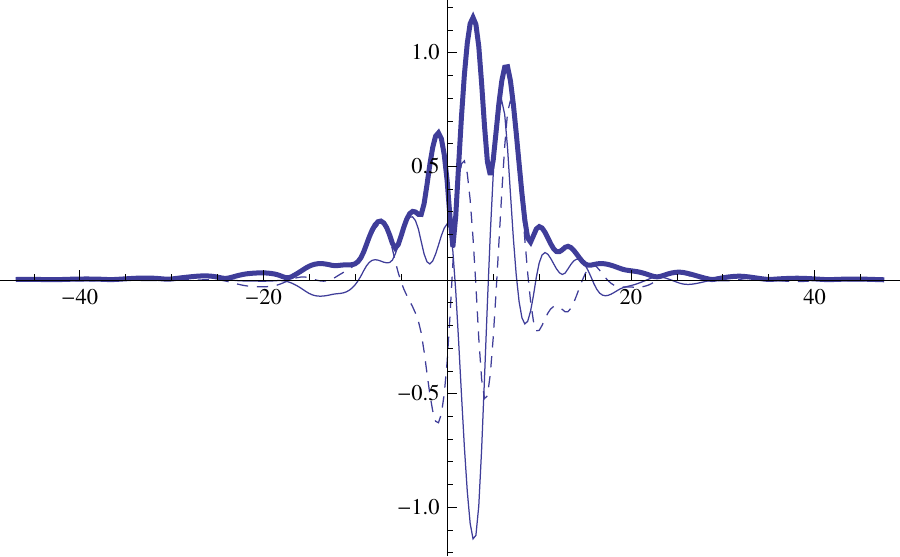}
\quad
\includegraphics[width = 0.45 \textwidth]{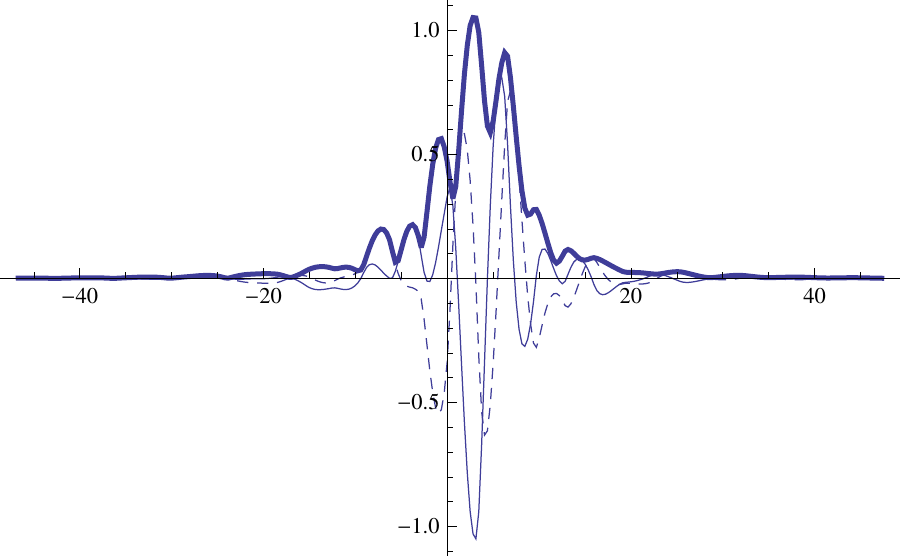}

(a) \quad \quad \quad \quad \quad\quad \quad \quad \quad \quad \quad \quad \quad \quad (b)

\bigskip

\includegraphics[width = 0.45 \textwidth]{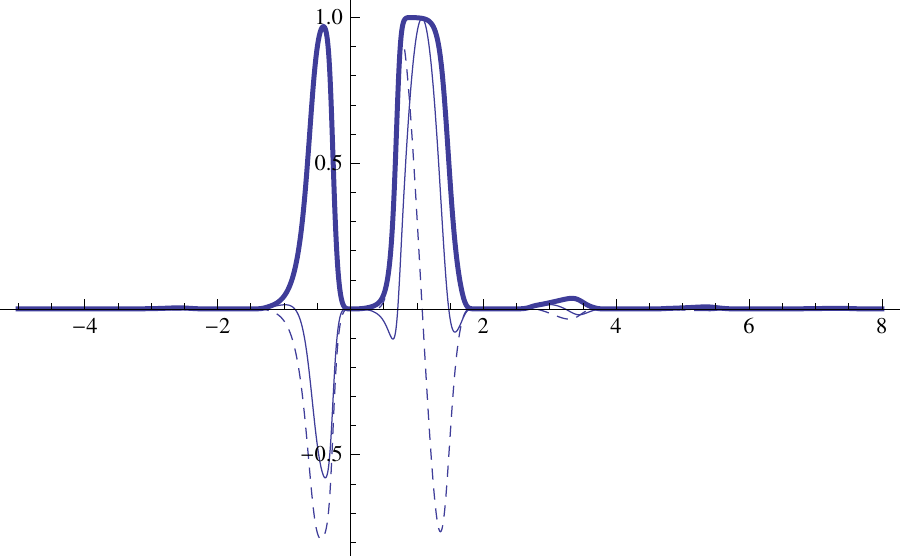}
\quad
\includegraphics[width = 0.45 \textwidth]{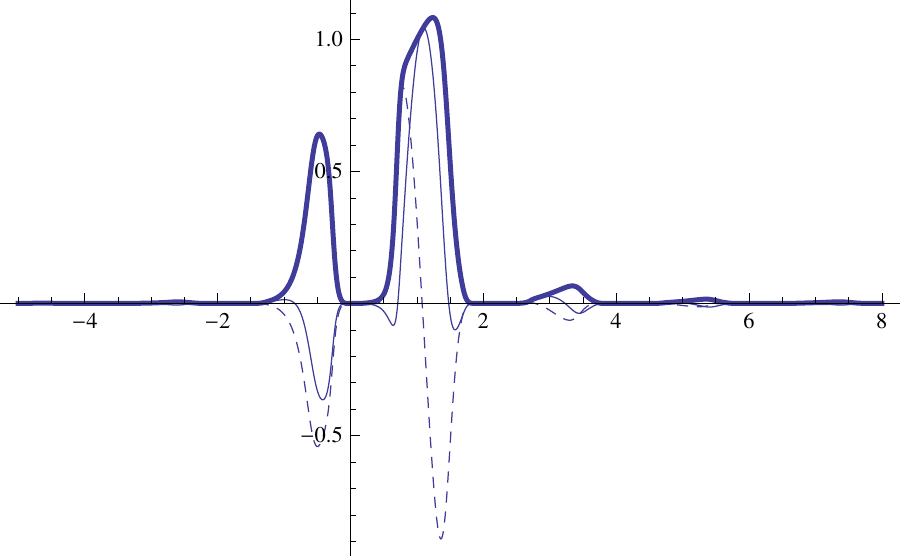}

(c) \quad \quad \quad \quad \quad\quad \quad \quad \quad \quad \quad \quad \quad \quad (d)

\caption{Time domain representation (a)  of the orthonormalized complex spline wavelet $\psi_{z,\perp}$ and (b) of its ridge variant $\Psi_{z,0,0}$.
Spectrum (c) of $\widehat{\psi}_{z,\perp}$ and (d) of  $\widehat{\Psi}_{z,0,0}$ for $z = 3.5 + i$. Thick lines indicate the modulus, thin lines the real part and dashed lines the imaginary part.}

\label{fig Spline ridge wavelets Frequenz}
\end{figure}

\begin{figure}[p]
\centering
\includegraphics[width = 0.45 \textwidth]{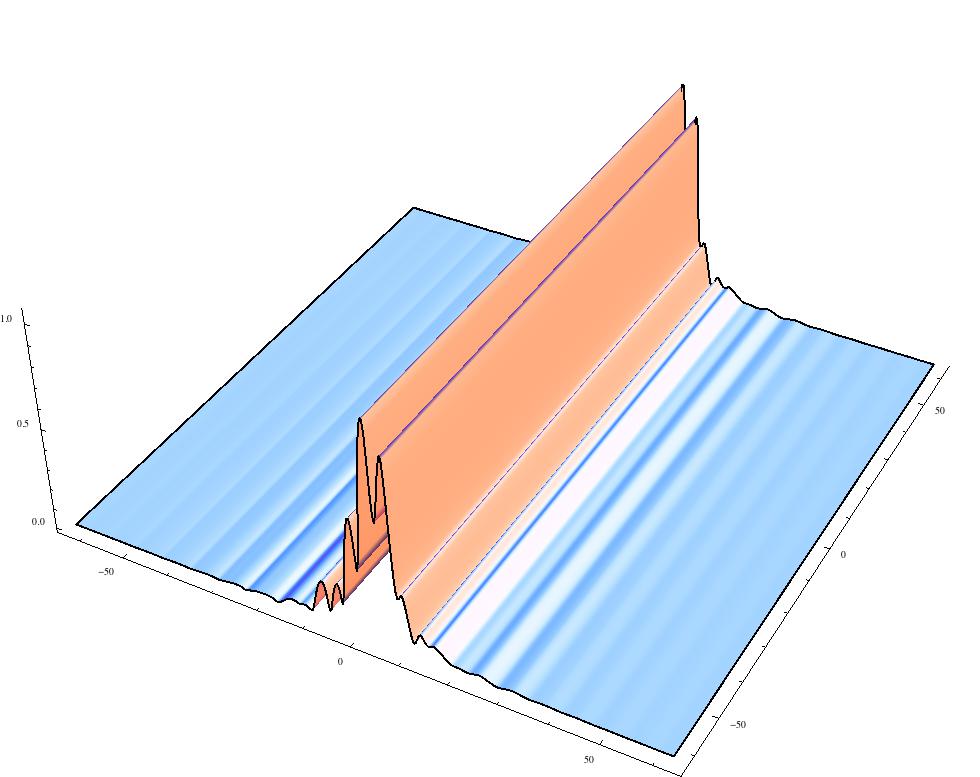}

 (a)

\bigskip

\includegraphics[width = 0.45 \textwidth]{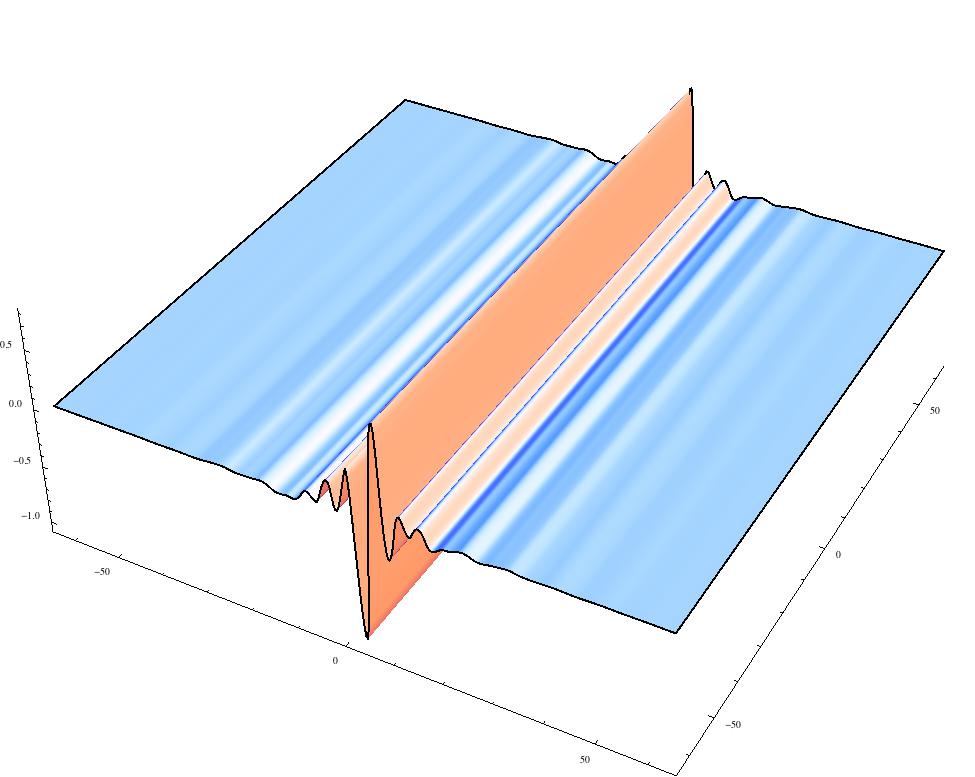}
\quad \quad
\includegraphics[width = 0.45 \textwidth]{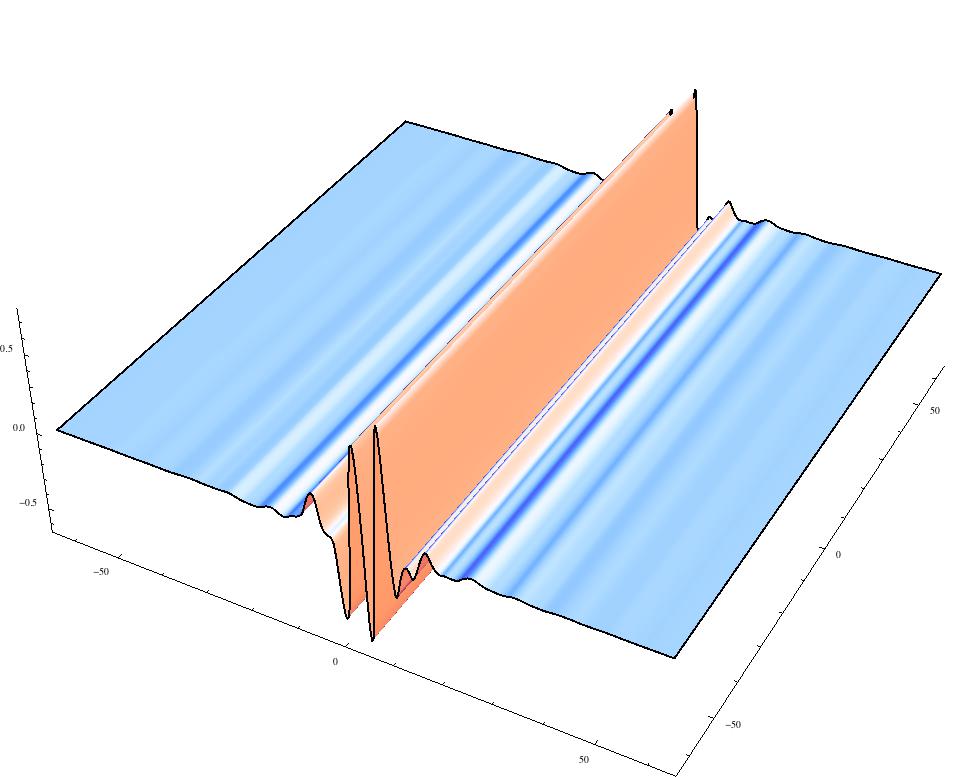}

(b) \quad \quad \quad \quad \quad \quad \quad \quad \quad \quad \quad \quad \quad \quad (c)

\bigskip

\caption{Orthonormal spline ridge wavelet $\Psi_{z,0,0} (u\cdot x)$ associated to the complex spline scaling function ${\beta}_{z,\perp}$ for $z=3.5 + i$ for $u = (1 \ 0)^{T}.$
 (a) Modulus, (b) real and (c) imaginary part. }

\label{fig Spline ridge wavelets}
\end{figure}

Although complex B-splines do not belong to the function class $\cas$, for the existence of integrals of the type \eqref{integral} that involve functions as defined in \eqref{126b}, it is not necessary that $g\in \cas$. It suffices that $f$ and $g$  belong to $L^{2}(\mr)$ or that additionally one of them belongs to a Sobolev space $H^{r}(\mr)$ with $r > \frac{n-1}{2}$, i.e., an appropriate subspace of $L^{2}(\mr)$. Choosing $z$ for $\beta_z$ so that $\re z > \frac{n}{2}$ these integrals exist. Therefore, we can define
\[
\Psi_{z; k, l, u} (x) := \Psi_{z; k, l} (u\cdot x) = \bDiff^{\frac{n-1}{2}} \psi_{z;k,l}(u\cdot x) = (|\cdot|^{\frac{n-1}{2}}\widehat{\psi_{z;k,l}})^\lor (u \cdot x),
\]
we can again apply Theorem \ref{126k} to obtain decompositions of  $L^{1}(\mathbb{R}^{n}) \cap L^{2}(\mathbb{R}^{n}).$ For illustrations, see Figures \ref{fig Spline ridge wavelets Frequenz} and \ref{fig Spline ridge wavelets}.

%a decomposition of $L^{2} \cap L^{1}(\mathbb{R}^{n})$.

%Alternatively, one could define dual wavelets by setting $\widehat{\widetilde{\psi}}_z := \widehat{\psi}_z/R_z$. Then,
%\[
%\Psi_{z; k, l, u} (x) := \Psi_{z; k, l} (u\cdot x) = \bDiff^{\frac{n-1}{2}} \psi_{z;k,l}(u\cdot x) = (|\cdot|^{\frac{n-1}{2}},\widehat{\psi_{z;k,l}})^\lor (u \cdot x)
%\]
%and
%$$
%\widetilde{\Psi}_{z; k, l, u} (x) = (|\cdot|^{\frac{n-1}{2}},\widehat{\widetilde{\psi}_{z;k,l}})^\lor (u \cdot x)
%$$
%generate dual continuous frames for $L^{2}(\mathbb{R}^{n})$.
\enx

\section{Discrete Representations}

In this section, we consider the cube,
\bes Q := [-1,1]^n \subset\mr^n,\ens
and functions $f\in L^2(Q)$.
We will present a discretization of the sphere $\bS^{n-1}$ which ultimately leads to a complete discrete representation of functions $f\in L^2(Q)$. This discretization was also considered in \cite{C} and is based on the concept of an $\eps$--net. It is one of several existing discretization methodologies (other choices include the methods in \cite{conradprestin}, \cite{ward}, and  \cite{brauchart}).
Let us recall the definition of a finite $\eps$--net.
\bd
Let $(X,d)$ be a metric space and a discrete set $N\subset X$. Given any $\epsilon >0,$ the set $N$ is called an $\eps$--net for $M$ if
\begin{itemize}
\item[(a)] $\inf\{d(y,y') \mid y\neq y'\in N\}\geq\eps$;
\item[(b)] $\inf\{r \mid \bigcup_{y\in N} \overline{B}_r(y)\supseteq X\} \leq\eps$, where $\overline{B}_r(y)$ denotes the closed ball of radius $r>0$ centered at $y$.
\end{itemize}
An $\eps$-net is called finite if $N$ is a finite set.
\ed
Note that since the sphere is compact, hence totally bounded, an $\eps$--net exists for $\bS^{n-1}$ for all $\eps > 0;$ see \cite{suther}.

We employ the following discretization procedure for $\bS^{n-1}$; see also \cite{C}.
\begin{itemize}
\item[(i)] Choose an $a_0 >1$, and discretize the scale parameter $a$ by the sequence $\{a_k := a_0^k \mid k \in I\}$, where $I := \{k\in \mz \mid k \geq k_0\}$ and $k_0\in \mz$ is selected appropriately.
\item[(ii)]  For $k\in I$, set $\eps_k := \frac12\,a_{k_0-k}$.
\item[(iii)] Let $S_k^{n-1}$ be an $\eps_k$-net of $\bS^{n-1}$ and require that following condition holds: There exist positive constants $c= c(n)$ and $C= C(n)$ so that for all $u\in \bS^{n-1}$ and all $\eps_k\leq r \leq 2$
\[
c \left(\frac{r}{\eps_k}\right)^{n-1} \leq \card (\{B_r(u) \cap \bS_k^{n-1}\})\leq C \left(\frac{r}{\eps_k}\right)^{n-1}.
\]

\end{itemize}
Note that for $0 < r \leq \eps_k$, $B_r(u) \subseteq B_{\eps_k} (u)$, and thus $\card \{B_r(u) \cap S_k^{n-1}\} \leq C$. It can be proved that the number of points $N_k$ in the $\eps_k$-net satisfies the following bounds:
\[
c \left(\frac{r}{\eps_k}\right)^{n-1} \leq N_k \leq C \left(\frac{r}{\eps_k}\right)^{n-1}.
\]
Next, we list the standing assumptions for this section.
\vskip 6pt\noindent
\textbf{General setup:} Let $g\in \cas$ and assume that
\begin{itemize}
\item[(i)] $\displaystyle{\inr \frac{|\widehat{g}(\gamma)|^2}{|\gamma|^n} d\gamma < \infty}$;
\item[(ii)]  $\displaystyle{\inf_{1 \leq |\gamma|\leq a_0} \,\sum_{k = 0}^\infty \,\left\vert \widehat{g}(a_0^{-k}\gamma)\right\vert^2 \,\left\vert a_0^{-k}\gamma\right\vert^{-2(n-1)} > 0}$;
\item[(iii)] $\left\vert \widehat{g}(\gamma)\right\vert \leq K |\gamma|^\alpha\,(1 + |\gamma|)^{-\beta}$, for some $K>0$, $\alpha > \frac{n-1}{2}$ and $\beta > \alpha + \frac{n+3}{2}$.
\end{itemize}
In particular, we remark that
\begin{itemize}
\item if condition (i) is satisfied, then  $G := \bDiff^{\frac{n-1}{2}} \,g$ satisfies the admissibility condition \eqref{ad1};
\item condition (ii) is satisfied if
\[
\displaystyle{\inf_{1 \leq |\gamma|\leq a_0} \,\sum_{k = 0}^\infty \,\left\vert \widehat{g}(a_0^{-k}\gamma)\right\vert^2 \,\left\vert a_0^{-k}\gamma\right\vert^{-(n-1)} > 0}.
\]
\end{itemize}

For the proof of our result  we need \cite[Theorem 4]{C}, which we state here in our notation:
\bt\cite[Theorem 4]{C}\label{defC}
Assume that the function $g\in C^1(\mr)$  satisfies the following two conditions:
\begin{itemize}
\item[$\bullet$] $\displaystyle{\inf_{1 \leq |\gamma|\leq a_0} \,\sum_{k = 0}^\infty \,\left\vert \widehat{g}(a_0^{-k}\gamma)\right\vert^2 \,\left\vert a_0^{-k}\gamma\right\vert^{-(n-1)} > 0}$;
\item[$\bullet$] $\left\vert \widehat{g}(\gamma)\right\vert \leq K |\gamma|^\alpha\,(1 + |\gamma|)^{-\beta}$, for some $K>0$, $\alpha > \frac{n-1}{2}$ and $\beta > 2 + \alpha$.
\end{itemize}
Then there exists $b_0 > 0$ so that for any $b < b_0$, we can find two constants $A, B > 0$ (depending on $g$, $a_0$, $b_0$, and $n$) so that, for any $f\in L^2 (Q)$,
\bee \label{802}
A \|f\|^2_{L^2(\mr)} \leq \sum_{k\in I}\,\sum_{u\in S_k^{n-1}}\,\sum_{\ell\in\mz} \vert \langle f, D_{a_k} T_{\ell b}\,G_{u}\rangle \vert^2 \leq B \|f\|^2_{L^2(\mr)}.
\ene
\et

We will now show that under the general setup and with the discretization of the unit sphere $\bS^{n-1}$ in term of the $\eps$--net introduced above, there exists a discrete frame for $L^2(Q)$.

\bt
\label{S 4.3}
Let $g\in \cas$ be as in the general setup and let $G := \bDiff^{\frac{n-1}{2}} \,g$. Then there exists a $b_0 > 0$ so that \eqref{802} holds for any given $b \in ]0, b_0],$ i.e., the orthogonal projection of
$\{D_{a_k} T_{\ell b}\,G_{u} \mid k\in I;\, \ell\in \mz;\,u\in S_k^{n-1}\}$ onto
$L^2(Q)$ forms a 
frame for $L^2(Q)$.
\et

\bp Let $G$ be defined as in \eqref{126b}. Then,
\begin{eqnarray*}
\left\vert \widehat{G}(a_0^{-k}\gamma)\right\vert^2 \,\left\vert a_0^{-k}\gamma\right\vert^{-2(n-1)} & = & \left\vert \widehat{g}(a_0^{-k}\gamma)\right\vert^2 \,\left\vert a_0^{-k}\gamma\right\vert^{n-1}\,\left\vert a_0^{-k}\gamma\right\vert^{-2(n-1)}\\
& = & \left\vert \widehat{g}(a_0^{-k}\gamma)\right\vert^2 \,\left\vert a_0^{-k}\gamma\right\vert^{-(n-1)},
\end{eqnarray*}
and, therefore,
\[
\inf_{1 \leq |\gamma|\leq a_0}\, \sum_{k= 0}^\infty \,\left\vert \widehat{G}(a_0^{-k}\gamma)\right\vert^2 \,\left\vert a_0^{-k}\gamma\right\vert^{-(n-1)} > 0.
\]
Furthermore,
\[
\vert \widehat{G}(\gamma)\vert = \left\vert \widehat{g}(\gamma)\right\vert \vert\gamma \vert^{(n-1)/2} \leq K |\gamma|^{\alpha + (n-1)/2}\,(1 + |\gamma|)^{-\beta}.
\]
Hence, the function $G$ satisfies the two conditions in Theorem \ref{defC} and the result follows.
\ep

\end{document}